\documentclass[preprint,12pt]{article}
\usepackage{amsmath,graphicx,color,algorithm,blkarray,authblk}

\begin{document}

\title{The elapsed time between two transient state observations for an absorbing Markov Chain}
\author[1]{Bianca~De~Sanctis}
\author[1,2]{A.P.~Jason~de~Koning}
\renewcommand\Affilfont{\itshape\small}

\affil[1]{Department of Biochemistry and Molecular Biology, University of Calgary}
\affil[2]{Alberta Children's Hospital Research Institute and Department of Medical Genetics, University of Calgary}
\maketitle

\begin{abstract}
Consider a system evolving according to an absorbing discrete-time Markov chain with known transition matrix. The state of the system is observed at two points in time, separated by an unknown number of generations. We are interested in calculating the expected elapsed time and its variance. We provide a novel, exact solution, which is computable from the fundamental matrix of a related absorbing Markov chain. This solution may be useful in population genetics for computing the expected age of a segregating allele without requiring diffusion approximation.
\end{abstract}

\section{Introduction}

Rigorously stated, the problem is as follows. Let $X(t)$ be an absorbing Markov chain with known transition matrix $P$. Let $S$ be the set of transient states. Assume $X(0)=s_i\in S$ and $X(T)=s_j \in S$ are known, where $T$ is a discrete random variable with positive integer support. What are $E(T | s_i,s_j,P)$ and $V(T | s_i,s_j,P )$? Intuitively, the elapsed time should be finite due to the guaranteed eventual absorption of the Markov chain. \\

Here we derive expressions for the expected value and variance of the elapsed time using the law of total expectation and the law of total variance, respectively. For higher moments, it is possible to apply the law of total cumulance. Ideally one could find the full conditional distribution of the random variable $T$, and we leave this as an open problem for future research. We note that this question only has an answer when the Markov chain is absorbing, and only when both observed states are transient. \\

The remainder of the paper is organized as follows. Sections \ref{expsec} and \ref{varsec} are derivations of the expected value and variance, respectively. Throughout these sections, we use values dependent on the transition matrix such as $H_{jj}$, the recurrence probability, $\tau_{ij}$, the mean first passage time conditional upon passage, and $v_{ij}$, the variance of the first passage time conditional upon passage. A derivation of these quantities is found in Section \ref{extras}. Section \ref{example} gives a simple numerical application of these results and Section \ref{discuss} gives an application to population genetics.

\section{Expected value}\label{expsec}

Here we derive $E(T | s_i,s_j,P)$, which from here on is referred to as $E(T)$. \\

Define a random variable $N$ to be the number of times that state $s_j$ has occurred while the observer has been gone and including the final observation. Notice that $N\geq 1$, which makes explicit the implicit assumption that at least one step has elapsed. Let $H_{jj}$ be the recurrence probability of transient state $s_j$, the probability that, starting in state $s_j$, the chain will ever reach $s_j$ again. Let $\tau_{ij}$ be the expected first conditional passage time from a transient state $s_i$ to a transient state $s_j$, i.e. the expected number of steps to reach a transient state $s_j$ from a transient state $s_i$, given that passage to $s_j$ occurs. The derivation of these values is described in Section \ref{extras}. \\

A special case of the law of total expectation allows us to condition as
\begin{align}\label{condition} E(T) = \sum_{n=1}^{\infty} E(T|N=n) P(N=n) \end{align}

First consider $P(N=n)$. Since $s_j$ was observed, we have $P(N\geq 1) = 1$. Now, the probability that $s_j$ was hit more than once is simply the probability that, given the chain was in $s_j$, the chain returns to $s_j$. Mathematically, $P(N>1) = H_{jj}$. Induction then gives $P(N>n) = H_{jj}^n$. Therefore, 

\begin{equation}\label{p} P(N=n) = P(N>n-1) - P(N>n) = H_{jj}^{n-1} - H_{jj}^{n} \end{equation}

Now we wish to find $E(T|N=n)$. The random variable $T|N=n$ can be partitioned into the first conditional passage time $t_{ij}$ that the process hits $s_j$ having started in $s_i$, and $n-1$ further conditional passage times $t_{jj}$ that the process returns to $s_j$ having already hit $s_j$. Then 

\begin{equation}\label{e} E(T|N=n) = E( t_{ij} + (n-1) t_{jj} ) = \tau_{ij} + (n-1)\tau_{jj}\end{equation}

by the linearity of expectation.\\

Plugging Equations \ref{p} and \ref{e} into Equation \ref{condition} gives

\begin{align}E(T) = \sum_{n=1}^{\infty} (\tau_{ij} + (n-1)\tau_{jj})(H_{jj}^{n-1} - H_{jj}^{n} ) \end{align}

Here we state some useful summation identities for $|x|<1$.
The second can be easily seen by differentiating the first, and likewise the third from the second.

\begin{align}\label{nnm1} \sum_{n=1} x^n = \frac{x}{1-x}, \; \; \; \sum_{n=1}^{\infty} nx^{n-1} = \frac{1}{(1-x)^2}, \; \; \; \sum_{n=1}^{\infty} n^2 x^{n} = \frac{x(x+1)}{(1-x)^3}\end{align}

A little algebra then yields

\begin{align} E(T) &= (\tau_{ij} - \tau_{jj} )\sum_{n=1}^{\infty} \left( H_{jj}^{n-1} - H_{jj}^{n} \right) + \tau_{jj} \sum_{n=1}^{\infty} n\left(H_{jj}^{n-1} - H_{jj}^n\right) \\&= \tau_{ij} + \tau_{jj} \left(\frac{H_{jj}}{1-H_{jj}} \right)\end{align}

This result is quite intuitive if we think about the number of returns to $s_j$ as a geometric distribution with failure probability $H_{jj}$.

\section{Variance}\label{varsec}
Again let $V(T|s_i,s_j,P) = V(T)$ for simplicity. Similar to the law of total expectation, a special case of the law of total variance states that 

\begin{align}\begin{split}\label{vari} V(T)  &= \sum_{n=1}^{\infty} V(T | N=n) P(N=n) + \sum_{n=1}^{\infty} E(T|N=n)^2 (1-P(N=n)) P(N=n) \\ &- 2 \sum_{n=1}^{\infty} \sum_{m=1}^{n-1} E(T | N=n) P(N=n) E(T | N=m) P(N=m) \end{split}\end{align}

We have all of these values from the last section except for $V(T|N=n)$. We can approach this similarly to $E(T | N=n)$. Partition $T$ into the first conditional passage time $t_{ij}$ and $n-1$ further conditional passage times $t_{jj}$ to get
\begin{align} V(T | N=n) = V(t_{ij} + (n-1) t_{jj})\end{align}

Because of independence due to the memoryless property of the Markov chain,
\begin{align}V(T|N=n) = V(t_{ij}) + (n-1)^2 V(t_{jj}) = v_{ij} + (n-1)^2 v_{jj} \end{align}

where $v_{ij}$ is the variance of the first passage time from $s_i$ to $s_j$, conditional on passage. The derivation of $v_{ij}$ is shown in Section \ref{extras}. \\

Plugging all relevant values into equation \ref{vari} yields

\begin{align}\begin{split}\label{fullvar} V(T) &= \sum_{n=1}^{\infty} (v_{ij} + (n-1)^2v_{jj})(H_{jj}^{n-1}-H_{jj}^{n}) \\& + \sum_{n=1}^{\infty} (\tau_{ij} + (n-1)\tau_{jj})^2(1-H_{jj}^{n-1}+H_{jj}^n)(H_{jj}^{n-1}-H_{jj}^{n}) \\&- 2\sum_{n=1}^{\infty} \sum_{m=1}^{n-1} (\tau_{ij} + (n-1)\tau_{jj}) (H_{jj}^{n-1}-H_{jj}^{n} ) (\tau_{ij} + (m-1)\tau_{jj})(H_{jj}^{m-1}-H_{jj}^{m}) \end{split} \end{align}

Making heavy use of the identities in Equation \ref{nnm1}, we obtain for the first sum in Equation \ref{fullvar}

\begin{align}\begin{split} \label{var1} \sum_{n=1}^{\infty} (v_{ij} + &  (n-1)^2v_{jj})(H_{jj}^{n-1} - H_{jj}^{n} )  = v_{ij} + v_{jj} + v_{jj} \left(\frac{H_{jj}+1}{(1-H_{jj})^2} \right) -\frac{2v_{jj}}{1-H_{jj}} \\
& = v_{ij} + v_{jj} \left( \frac{3H_{jj}-1}{(H_{jj}-1)^2} + 1 \right)
\end{split} \end{align}

For the second sum we get

\begin{align}\begin{split} \label{var2} \sum_{n=1}^{\infty} (\tau_{ij} & + (n-1)\tau_{jj})^2(1-H_{jj}^{n-1}+H_{jj}^n)(H_{jj}^{n-1}-H_{jj}^{n}) \\
= &(\tau_{ij}^2-2\tau_{ij}\tau_{jj} + \tau_{jj}^2)\left( \frac{2H_{jj}}{H_{jj}+1} \right) + \left(2\tau_{ij}\tau_{jj} - 2\tau_{ij}^2 \right) \left( \frac{-H_{jj}^2-3H_{jj}}{(H_{jj}-1)(H_{jj}+1)^2} \right) \\ &+ \tau_{jj}^2 \left( \frac{H_{jj}^4 + 5H_{jj}^3 + 5H_{jj}^2 + 5H_{jj}}{(H_{jj}-1)^2(H_{jj}+1)^3} \right)
\end{split} \end{align}

For the third sum we get
\begin{multline} \begin{split} \label{varr3} - 2\sum_{n=1}^{\infty} & \sum_{m=1}^{n-1} (\tau_{ij} + (n-1)\tau_{jj}) (H_{jj}^{n-1}-H_{jj}^{n} ) (\tau_{ij} + (m-1)\tau_{jj})(H_{jj}^{m-1}-H_{jj}^{m}) \\
&= -2H_{jj}  (\tau_{ij}^2 H_{jj}^4 - 2\tau_{ij}^2 H_{jj} ^2 + \tau_{ij}^2 - 2\tau_{ij}\tau_{jj} H_{jj} ^4 - \tau_{ij}\tau_{jj}H_{jj} ^3 + 3\tau_{ij}\tau_{jj} H_{jj} ^2 
\\ &
+ \tau_{ij} \tau_{jj} H_{jj}  - \tau_{ij} \tau_{jj} + \tau_{jj}^2 H_{jj} ^4 + \tau_{jj}^2 H_{jj} ^3  - 2\tau_{jj}^2 H_{jj} ^2 - 2\tau_{jj}^2 H_{jj}  - 2\tau_{jj}^2 )\\ 
& \big/ ((H_{jj} -1)^2(H_{jj} +1)^3)
\end{split} \end{multline}

Thus we have computed the variance as the sum of Equations \ref{var1}, \ref{var2} and \ref{varr3}. 

\section{Derivation of $H_{jj}$, $\tau_{ij}$ and $v_{ij}$}\label{extras}

We begin with the goal of deriving $\tau_{ij}$ and $v_{ij}$, and will derive $H_{jj}$ along the way. Since we are in an absorbing Markov chain, we cannot use the mean first passage time commonly denoted $\mu_{ij}$. This value is dependent upon the first passage occurring, which is not guaranteed since our chain is absorbing. Indeed, if we attempted to compute it, we would get $\infty$. This is not the value we're looking for. \\

In our case, we need something we will call the conditional mean first passage time, or $E(t_{ij}) = \tau_{ij}$. We will call the conditional first passage time variance $V(t_{ij} ) = v_{ij}$. The first is the expected time until the chain hits the state $s_j$, conditional on starting in state $s_i$ and on hitting $s_j$. We can compute this and the variance $v_{ij}$ using known procedures when $i\ne j$, and once we have done this, the derivation is simple for the case when $i=j$. \\

First we will consider the case when $i \ne j$. Since, when considering first passage times, it does not matter what happens after we reach state $s_j$, the probabilities $p_{j \cdot} $ out of state $s_j$ are irrelevant. This justifies temporarily redefining $s_j$ as an absorbing state in order to perform our calculations. Once we do this, we are simply looking for the conditional expected time to absorption and the variance of the conditional time to absorption, for which equations are known. \\

The procedure to compute $\tau_{ij}$ and $v_{ij}$ where $i\ne j$ is thus as follows, and is more or less taken from Kemeny and Snell's book \cite{markovbook}. First, replace the $j$-th row of $P$ with a $1$ in the $j$-th coordinate and zeros everywhere else, so $s_j$ becomes an absorbing state. Let $t$ be the number of transient states of this modified Markov chain, one less than the number of transient states of $P$, and let $r$ be the number of absorbing states of this modified Markov chain, one more than the number of absorbing states of $P$. Now, reorganize the modified Markov chain transition matrix as
\begin{align} \label{split} \tilde{P} = \left(\begin{array}{cc} \tilde{Q} & \tilde{R} \\ 0 & I_r \end{array}\right)\end{align}
where $\tilde{Q}$ is a $t\times t$ matrix of transient-to-transient state transitions, $\tilde{R}$ is a nonzero $t\times r $ matrix of transient-to-absorbing state transitions, $I_r$ is the $r\times r$ identity matrix of absorbing-to-absorbing state transitions, and $0$ is an $r \times t$ matrix of absorbing-to-transient state transitions where every entry is $0$. Then
\begin{equation}\label{fundamental} \tilde{N} = \sum_{k=0}^{\infty}\tilde{Q}^k = (I_t-\tilde{Q})^{-1}\end{equation}
is called the fundamental matrix of the modified Markov chain. The probability of absorbing in state $l$ having started in state $i$ is the $(i,l)$th entry of the matrix 
\begin{equation}\label{B} \tilde{B} = \tilde{N} \tilde{R} \end{equation}
where $\tilde{R}$ is defined as in Equation \ref{split}. In our context, the matrix $\tilde{B}$ contains, in the $i$th row and the column corresponding to the new absorbing state $s_j$, the probability that state $s_j$ is ever reached, conditional on starting in state $s_i$. In other words, the column corresponding to the absorbing state $s_j$ contains the hitting probabilities $H_{ij}$ for $i \ne j$. We can then obtain $H_{jj}$ by conditioning on the first step out of $s_j$ as
\begin{align}\label{hjj} H_{jj} = p_{jj} + \sum_{k \in S, k \ne j} p_{jk}H_{kj} \end{align}

Next define $\tilde{D}_j$ as a diagonal matrix with diagonal entries defined as $\tilde{b}_{ij}$, $i \in S$, for the fixed absorbing state $s_j$. Then the expected number of steps before absorption having started in transient state $i$ and conditional on absorbing in state $j$ is the $i$th entry of
\begin{equation}\label{cond2} \mathbf{\tilde{t}} = \tilde{D}_j^{-1} \tilde{N} \tilde{D}_{j} \mathbf{1} \end{equation}
where $\mathbf{1}$ again is a vector of $1$s. In other words, the $i$-th entry of $\mathbf{\tilde{t}}$ is $\tau_{ij}=E(t_{ij})$ for $i \ne j$. The variance on this is the $i$th entry of
\begin{equation}\label{variancee} \mathbf{\tilde{v}} = (2\tilde{D}_j^{-1} \tilde{N} \tilde{D}_{j}-I_{t})\mathbf{\tilde{t}} - \mathbf{\tilde{t}}_{sq}\end{equation} 
where $\mathbf{\tilde{t}}_{sq}$ is $\mathbf{\tilde{t}}$ with every entry squared. In other words, the $i$th entry of $\mathbf{\tilde{v}}$ is $v_{ij}$ for $i \ne j$. Note that for a fixed $j$, we need to compute the entirety of the vectors $\mathbf{\tilde{t}}$ and $\mathbf{\tilde{v}}$ to use in the derivation of $\tau_{jj}$ and $v_{jj}$ below.\\

Now we show how to obtain $E(t_{ij}) = \tau_{ij}$ when $i=j$, which from here on is referred to as $\tau_{jj}$. This is the mean first expected recurrence time from $s_j$ to $s_j$, given recurrence. Of course this is not obtainable from the above procedure, since the above procedure requires defining $s_j$ as an absorbing state. If the chain is in state $s_j$, it will either remain there with probability $p_{jj}$ (with this probability taken from the original transition matrix $P$), in which case it will have taken one step, or it will go to a different state with probability $p_{ji}$, $i\ne j$, in which case it will take an expected $\tau_{ij}$ steps to return to $j$, plus the one step it took to get away. This is almost correct, but notice that the probabilities $p_{jj}$ and $\sum_{i \in S, i \ne j} p_{ji}$ will not necessarily add up to $1$. This is because there are more absorbing states to consider. However, we want $\tau_{jj}$ to be conditional on eventually hitting $s_j$, so we must condition on not hitting the original absorbing states, that is we must modify our probabilities so that they add up to $1$. This can be done by dividing the equation by one minus the probability that the process absorbs in the original absorbing states of $P$, i.e. those $s_i$ such that $i\not \in S$.
So we have

\begin{equation} \tau_{jj} = E(t_{jj}) = E \left( \frac{ p_{jj} \cdot 1 + \sum_{i \in S, i \ne j} p_{ji} (t_{ij} + 1) } {1 - \sum_{i \not \in S} p_{ji}} \right) = \frac{ p_{jj}+ \sum_{i \in S, i \ne j} p_{ji} (\tau_{ij}+1)} {1- \sum_{i \not \in S} p_{ji}} \end{equation}

where $\tau_{ij}$, $i\ne j$, are taken from Equation \ref{cond2}. Likewise, we can compute variance as

\begin{equation} v_{jj} = V(t_{jj} ) = V \left( \frac{ p_{jj} \cdot 1 + \sum_{i \in S, i \ne j} p_{ji} (t_{ij} + 1) } {1 - \sum_{i \not \in S} p_{ji}} \right) = \frac{ \sum_{i\in S, i \ne j} p_{ji}^2 v_{ij} }{\left( 1- \sum_{i \not \in S} p_{ij} \right) ^2 }
\end{equation}

where $v_{ij}$, $i \ne j$, are taken from Equation \ref{variancee}.

\section{Example}\label{example}
Consider the absorbing Markov chain with transient states $S=\{s_1,s_2\}$ and absorbing states $\{s_0,s_3\}$ with transition matrix
\begin{align*} P = \left(\begin{array}{cccc} 1 & 0 & 0 & 0 \\ 0.5 & 0 & 0.5 & 0 \\ 0 & 0.5 & 0 & 0.5 \\ 0 & 0 & 0 & 1 \end{array}\right)\end{align*}
Let's assume that $i=1$, $j=2$, i.e. the first observation was $s_1$ and the later observation was $s_2$. We temporarily define $s_2$ as an absorbing state and rearrange the modified matrix according to Equation \ref{split}. The columns and rows of the modified matrix are now in the order $s_1,s_0,s_2,s_3$. 
\begin{align*} \tilde{P} = \left(\begin{array}{cccc} 0 & 0.5 & 0.5 & 0 \\
0 & 1 & 0 & 0 \\
 0 & 0 & 1 & 0 \\
 0 & 0 & 0 & 1 \end{array}\right)\end{align*}
We have $\tilde{Q} = (0)$ and $\tilde{R} = (0.5, 0.5, 0)$. Then
\begin{align*} \tilde{N} &= (I-\tilde{Q})^{-1} = (1)^{-1} = (1) \\
\tilde{B} &= (1) (0.5, 0.5, 0) = (0.5, 0.5, 0) = (H_{10}, H_{12}, H_{13} ) \\
H_{22} &= p_{22} + p_{21}H_{12} = 0 + 0.5 \times 0.5 = 0.25 \\
\tilde{D}_2 &= (0.5) = (H_{12})  \\
\tilde{\mathbf{t}} &= \tilde{D}_2^{-1} \tilde{N} \tilde{D}_{2} \mathbf{1} = (0.5)^{-1} (1) (0.5) (1) = (1) = (\tau_{12}) \\
\tilde{\mathbf{v}} & = (2(0.5)^{-1} (1) (0.5) - (1) ) (1) - (1) = (0) = (v_{12}) \\
t_{22} &= \frac{ p_{22} + p_{21} (\tau_{12} + 1) }{1- p_{23} - p_{24} } = \frac{0 + 0.5 (1 + 1)} {1 - 0 - 0.5 - 0} = 2 \\
v_{22} &= \frac{ p_{21}^2 v_{12} } {(1- p_{23} - p_{24})^2 } = \frac{0.5^2 \times 0} {1 - 0.5 - 0} = 0
 \end{align*}
All of these results are quite intuitive. From the original transition matrix $P$, it is clear that, having started in state $s_1$, the direct path is the only possible path to $s_2$. It has length $\tau_{12} = 1$, and since there are no other possibilities, we have $v_{12} = 0$. The hitting probability $H_{12}$ also makes sense, because the only possible path has probability $p_{12} = 0.5$. Likewise, the only possible path from $s_2$ back to $s_2$ is through $s_1$. This has a length of $\tau_{22} = 2$, and since there are no other possibilities, we again have $v_{22} = 0$. Also, since there are no other possible paths, the hitting probability $H_{22}$ has the probability of this path $p_{12}p_{21} = 0.5^2 = 0.25$. \\

Using these, we can compute $E(T)$ and $V(T)$. We do not show the cumbersome calculations for $V(T)$. 
\begin{align*} E(T) &= \tau_{12} + \tau_{22} \left( \frac{H_{22}} {1-H_{22}} \right) = 1 + 2 \left( \frac{.25}{1-.25} \right) = \frac{5}{3} \\
V(T) &= \frac{4744}{375} =  12.65066666... 
 \end{align*}

The variance in this example is  large compared to the expected value and we believe this is likely to often be the case, since there are many possible endpoint-conditioned trajectories. Despite that our solution is exact, intuitively, one pair of transient state observations may not provide much information about the elapsed time. Future research might therefore productively extend these results by considering multiple observations of the same realization, or by considering observations from multiple realizations of the same underlying process.

\section{Application to Population Genetics}\label{discuss}

The number of individuals that carry a mutant allele within a population can be described by an absorbing Markov chain such as the Wright-Fisher model \cite{fisher,wright}, where each state corresponds to the frequency of the allele in a diploid population. The transition matrix can be described completely given the population size, selection and dominance coefficients, and mutation rate to and from the allele, all of which can be estimated from data. The age of an allele can be calculated within our framework by assuming that it originates as a single mutation within one individual, giving $X(0)=s_i$ where $i=1$. The frequency of the allele of interest in the current population is $j$, giving $X(T)=s_j$. Using the Wright Fisher model to construct the transition matrix, our method allows calculation of the expected number of generations that the allele has existed, or the expected age of the allele. \\

The problem of calculating the expected age of an allele has been extensively studied for over 40 years. Many approximate results exist under different mathematical frameworks. The first of these results was obtained by Kimura and Ohta in 1973 \cite{kimohta} using a technique called diffusion theory, a continuous-time approximation to a Markov chain. Since then, more general approximations have been published using diffusion theory and other methods which include mutation, migration and selection \cite{review,slatbayes,malas}. However, we believe our approach is the first to give a method to calculate the expected allele age exactly. By virtue of this approach being completely general, it also uniquely allows for arbitrary selection, mutation, dominance, and population size, as these can all easily be included in the model without an increase in complexity. A potential shortcoming, however, is that it is not obvious how non-equilibrium demography could be accounted for. We leave this area for future research.

\section*{Acknowledgements}
This work was partly supported by a Discovery Grant  from the National Sciences and Engineering Research Council of Canada
(NSERC),  and by an NSERC\ USRA award\ (BDS). The authors also gratefully acknowledge infrastructure support from the Canada Foundation for Innovation and the Alberta\ Children's Hospital Research Institute.

\end{document}